\documentclass[preprint]{imsart}
\usepackage[official]{eurosym}
\usepackage{amssymb,amsmath}
\usepackage{epic,eepic}
\usepackage{array}
\usepackage{amsthm}
\usepackage{graphicx}
\usepackage[top=3cm, bottom=5cm, left=3cm, right=3cm]{geometry}
\usepackage{imsart}
\usepackage{fullpage,amsmath}
\usepackage{amssymb}
\usepackage{bbding}
\usepackage{graphicx}

\newtheorem{theorem}{Theorem}
\newtheorem{lemma}[theorem]{Lemma}
\newtheorem{corollary}[theorem]{Corollary}
\newtheorem{definition}[theorem]{Definition}
\newtheorem{remark}[theorem]{Remark}

\numberwithin{theorem}{section} \numberwithin{equation}{section}

\begin{document}

\title{Bernstein-type approximation of set-valued functions in the symmetric difference
metric\thanksref{T1}}
 \runtitle{Bernstein-type approximation of
set-valued functions in the symmetric difference metric}
\author{\fnms{Shay} \snm{Kels}\thanksref{t1}\ead[label=e1]{kelsisha@tau.ac.il}}
\and
\author{\fnms{Nira} \snm{Dyn}\ead[label=e2]{niradyn@tau.ac.il}}
\address{
School of Mathematical Sciences, Tel-Aviv University, Tel-Aviv,
Israel\\
\printead{e1,e2}} \runauthor{Shay Kels and Nira Dyn}
\thankstext{T1}{This is a preprint version of the paper.}
\thankstext{t1}{corresponding author}
\date{}
\maketitle

\begin{abstract}

We study the approximation of univariate and multivariate set-valued
functions (SVFs) by the adaptation to SVFs of positive sample-based
approximation operators for real-valued functions. To this end, we
introduce a new weighted average of several sets and study its
properties. The approximation results are obtained in the space of
Lebesgue measurable sets with the symmetric difference metric.

In particular, we apply the new average of sets to adapt to SVFs the
classical Bernstein approximation operators, and show that these
operators approximate continuous SVFs. The rate of approximation of
H{\"o}lder continuous SVFs by the adapted Bernstein operators is
studied and shown to be asymptotically equal to the one for
real-valued functions. Finally, the results obtained in the metric
space of sets are generalized to metric spaces endowed with an
average satisfying certain properties.

\end{abstract}
\maketitle
\section{Introduction}
Set-valued functions (SVFs) have various applications in
optimization, control theory, mathematical economics and other
areas. The approximation of SVFs from a finite number of samples has
been the subject of several recent research works
(\cite{baier2011set},\cite{dyn2007approximations},\cite{dyn2006approximations},\cite{kels2011subdivision})
and reviews (\cite{dynapproximation},\cite{muresanset}).

In order to adapt to SVFs sample-based approximation methods known
for real-valued functions, it is required to define linear
combinations of two or more sets. For most approximation methods it
is sufficient to consider linear combinations with weights summing
up to one, while for positive approximation operators only convex
combinations (non-negative weights summing up to one) are
considered. We term convex linear combinations as weighted averages.

In case of data sampled from a SVF mapping real-numbers to convex
sets, methods based on the classical \emph{Minkowski sum} of sets
can be used for the approximation
\cite{dyn2000spline,vitale1979approximation}. In this approach, sums
of numbers in positive operators for real-valued approximation are
replaced by Minkowski sums of sets. A generalization to sets which
are either convex or differences of convex sets is done in
\cite{baier2001differences}, where convex sets are embedded into the
Banach space of \emph{directed sets}. This approach allows to apply
existing methods for the approximation in Banach spaces
\cite{baier2011set}.

Approximation of set-valued functions mapping real-numbers to
general sets is a more challenging task. In this case, methods based
on Minkowski sum of sets fail to approximate the sampled function
\cite{vitale1979approximation,dyn2005set}, and other weighted
averages of sets are needed.

Artstein \cite{artstein1989pla} introduced a weighted average of two
sets with the property that the Hausdorff metric between the average
and any of the averaged sets changes linearly with the weight of the
average. This average was later termed as the \emph{metric average}
of sets. An extension of the metric average to the weighted average
of several sets, named the \emph{metric linear combination}, is
given in\cite{dyn2007approximations}. The metric linear combinations
were used in \cite{dyn2007approximations} to adapt to sets positive
and non-positive approximation operators, with the approximation
error measured in the Hausdorff metric. However, the metric linear
combination is applicable only to ordered sequences of sets, which
limits its usage to the approximation of univariate SVFs.

As it is noticed in \cite{artstein1989pla}, the particular choice of
a metric is crucial to the construction and analysis of set-valued
approximation methods. While previous works develop and analyze
set-valued approximation methods in the metric space of compact sets
endowed with the Hausdorff metric, we consider here the
approximation problem in the metric space of Lebesgue measurable
sets with the symmetric difference metric\footnote{The measure of
the symmetric difference is only a pseudo-metric on Lebesgue
measurable sets. The metric space is obtained in a standard way as
described in Section \ref{sectionPreliminaries}}. The symmetric
difference metric allows to obtain approximation results for a wider
class of functions, as is demonstrated in
\cite{kels2011subdivision},  where set-valued subdivision techniques
are investigated.

In this work, we consider the adaptation of positive sample-based
approximation operators  to univariate and multivariate SVFs . The
adaptation is based on a new weighted average of several sets,
termed the \emph{partition average}, which is studied in details.

As it is very well known, the concept of the weighted average of
numbers is closely related to that of the \emph{mathematical
expectation} of a \emph{ discrete random variable}. Similarly, a
weighted average of several sets may be interpreted as the
expectation of a \emph{random set} \cite{molchanov2005theory}. We
use tools from the theory of random sets to prove properties of the
partition average of sets.

First, we adapt to SVFs the classical Bernstein operators, and show
that these operators approximate continuous SVFs. Furthermore,  we
consider the rate of approximation of H{\"o}lder continuous SVFs by
set-valued Bernstein operators, and obtain a result for SVFs
analogous to that of Kac
\cite{kac1938remarque,kac1939reconnaissance} for real-valued
functions. Moreover, we show that the adaptation to SVFs of the
classical de Casteljau's algorithm (see, e.g.
\cite{farin2002curves}, Chapter 4) yields another sequence of
adapted operators having the same rate of approximation as that for
the adapted Bernstein operators.

The results for Bernstein operators are then extended to general
positive sample-based operators. Moreover, we study the application
of positive sample-based operators to monotone SVFs, and show that
the adapted operator is monotonicity preserving if and only if the
corresponding operator for real-valued functions is monotonicity
preserving.

Due to the commutativity of the partition average of sets, the
results are easily generalized to approximation operators for
multivariate SVFs. Finally, we generalize the approximation results
to functions with values in metric spaces endowed with a weighted
average, with properties similar to those of the partition average
of sets.

The structure of this work is as follows. In Section
\ref{sectionPreliminaries}, we survey definitions and results
relevant to our work. In Section \ref{sectionPartitionAverage}, we
introduce the partition average of sets and study its properties. In
Section \ref{sectionBernstein}, we adapt to sets the Bernstein
approximation operators. In section \ref{sectionCastelio}, we study
another type of set-valued Bernstein operators, obtained by adapting
to sets of the de Casteljau's algorithm . In Section
\ref{sectionOperators}, we consider the adaptation to SVFs of
positive sample-based operators. In Section \ref{sectionMonotone},
we discuss the approximation of monotone SVFs. The approximation of
multi-variate SVFs is the subject of Section \ref{sectionMulti}.
Finally in Section \ref{sectionGeneral}, we generalize the results
to functions with values in general metric spaces.

\section{Preliminaries}\label{sectionPreliminaries}
\subsection{Sets and the symmetric difference metric}
We denote by $\mu$ the $m$-dimensional \emph{Lebesgue measure} and
by $\mathfrak{L}$ the collection of \emph{Lebesgue measurable
subsets} of $\mathbb{R}^m$ having finite measure. The \emph{set
difference} of two sets $A,B$ is
\begin{equation*}
A \setminus B = \left\{ {p:p \in A,p \notin B} \right\} \ ,
\end{equation*}
and the \emph{symmetric difference} is defined by
\begin{equation*}
A\Delta B = A\setminus B \bigcup B\setminus A \ .
\end{equation*}
 The \emph{measure of the symmetric difference} of $A,B \in \mathfrak{L}$,
\begin{equation*}
d_\mu  \left( {A,B} \right) = \mu \left( {A\Delta B} \right) \ ,
\end{equation*}
induces a pseudo-metric on $\mathfrak{L}$, and  $\left(
{\mathfrak{L},d_\mu  } \right)$ is a complete metric space by
regarding any two sets $A, B$ such that $\mu \left( {A\Delta B}
\right) = 0$ as equal (\cite{halmos1974measure}, Chapter 8). For
$A,B \in \mathfrak{L}$, such that $B \subseteq A$, it is easy to
observe that
\begin{equation}\label{IncludedDistance}
d_\mu  \left( {A,B} \right) = \mu \left( {A\setminus B} \right) =
\mu \left( A \right) - \mu \left( B \right) \ .
\end{equation}

We use the notation ${\operatorname{ci}}\left( A \right)$ for the
\emph{closure} of the \emph{interior} of $A$. A bounded set $A$,
such that $A = \operatorname{ci}\left(A\right)$ is called
\emph{regular compact}. Regular compact sets are closed under finite
unions, but not under finite intersections, yet for $A,B$ regular
compact sets such that $B \subset A$,
\begin{equation}\label{ciInclusion}
A \bigcap B = B = \operatorname{ci}\left( {A \bigcap B} \right) \ .
\end{equation}

We recall that a set $A \in \mathfrak{L}$ is \emph{Jordan
measurable} if and only if its \emph{boundary} has zero Lebesgue
measure. Jordan measurable sets are denoted by $\mathbb{J}$. We
recall that $\mathbb{J}$ is closed under finite unions and finite
intersections. Note that for $A \in \mathbb{J}$,
\begin{equation}\label{jordanCIA}
\mu \left( A \right) = \mu \left( {{\operatorname{ci}}\left( A
\right)} \right) \ .
\end{equation}
Moreover for $B_0 ,...,B_n  \in \mathbb{J}$,
\begin{equation}\label{ciOfUnion}
{\text{ci}}\left( {\bigcup\limits_{i = o}^n {B_i } } \right) =
\bigcup\limits_{i = 0}^n {{\text{ci}}\left( {B_i } \right)} \ .
\end{equation}
We denote by $\mathfrak{J}$ the subset of $\mathbb{J}$ consisting of
\emph{regular compact sets}. Notice that for any $A,B \in
\mathfrak{J}$, $d_\mu  \left( {A,B} \right) = 0$ implies $A = B$,
therefore $d_\mu$ is a metric on $\mathfrak{J}$. In particular, the
empty set $\phi$ is in $\mathfrak{J}$, and it is the only set in
$\mathfrak{J}$ having zero measure. Note that by its definition
$\mathfrak{J}$ is closed under finite unions.

\subsection{Real-valued Bernstein approximation}\label{SubSectRealBernstein}
For a function $f:\left[{0,1} \right] \to \mathbb{R}$, the
\emph{Bernstein polynomial} of degree $n$ is
\begin{equation}\label{realBernstein}
B_n \left( {f,x} \right) = \sum\limits_{i = 0}^n {\left(
{\begin{array}{*{20}c}
   n  \\
   i  \\
 \end{array} } \right)} x^i \left( {1 - x} \right)^{n - i} f\left( {\frac{i}
{n}} \right) \ .
\end{equation}
The mapping $f\left(\cdot\right) \to B_n \left( {f,\cdot} \right) $
is called the \emph{Bernstein operator}. An extensive exposition of
Bernstein polynomials is given in \cite{devore1993constructive}.

Obviously one can interpret (\ref{realBernstein}) as the weighted
arithmetic average of the values $f\left( {\frac{j} {n}} \right)$.
The probabilistic nature of the Bernstein polynomials is also well
known. It can be recognized by interpreting the weights,
\begin{equation}\label{Bn}
b\left( {n,x;i} \right) = \left( {\begin{array}{*{20}c}
   n  \\
   i  \\

\end{array} } \right)x^i \left( {1 - x} \right)^{n - i} \ ,
\end{equation}
 as point  probabilities of a binomial distribution with parameters $n$ and
 $x$.

The polynomials $B_n \left( {f, \cdot } \right)$ are the basis of
Bernstein's proof of the Weierstrass Approximation Theorem
(\cite{bernstein1912demonstration}, see
\cite{levasseur1984probabilistic} for a modern presentation). Using
Bernstein polynomials the theorem can be formulated as
\begin{theorem}\label{Weierstrass} Let $f:\left[ {0,1} \right] \to \mathbb{R}$ be a
continuous function, then for any $\varepsilon > 0$ there exists $N
> 0$, such that for all $n \geq N$ and all $x \in \left[ {0,1} \right]$ ,
\begin{equation*}
\left| {f\left( x \right) - B_n \left( {f,x} \right)} \right| <
\varepsilon {\text{ }} \ .
\end{equation*}
\end{theorem}

A stronger version of the above theorem for H{\"o}lder continuous
functions is due to Mark Kac
(\cite{kac1938remarque,kac1939reconnaissance}, see
\cite{mathe1999approximation} for a modern presentation). We denote
by $\operatorname{Lip}\left(L,\nu\right)$ the class of H{\"o}lder
continuous functions with exponent $\nu$ and constant $L$, defined
on $[0,1]$, namely functions satisfying,
\begin{equation}\label{HolderFunctionEq}
\left| {f\left( x \right) - f\left( y \right)} \right| \leq L\left|
{x - y} \right|^\nu \ , x,y \in [0,1] \ .
\end{equation}

\begin{theorem}\label{Kac}Let $f \in \operatorname{Lip}\left(L,\nu\right)$,
then
\begin{equation*}
\left| {f\left( x \right) - B_n \left( {f,x} \right)} \right| \leq
L\left( {\frac{{x\left( {1 - x} \right)}} {n}} \right)^{\nu /2} \ .
\end{equation*}
\end{theorem}

Our adaptation of Bernstein operators to SVFs is based on the new
average of sets introduced in Section \ref{sectionPartitionAverage}.
To obtain the relevant properties of the new average of sets, we
give it a probabilistic interpretation using the notion of a
\emph{random closed set}, which is discussed together with basic
relevant results in the next subsection.

\subsection{Random sets}
We proceed with a few definitions regarding \emph{random sets}. The
following definitions and results are adapted from
\cite{molchanov2005theory}, which provides a thorough account of
random sets theory.

Here we denote by $\mathcal{F}$ the collection of closed subsets of
$\mathbb{R}^m$.
\begin{definition}Let $\left\{ {\Omega ,\mathfrak{F},\Pr} \right\}$
be a probability space. A map  $X:\Omega \to \mathcal{F}$ is called
a random closed set, if for every compact set $K \subset
\mathbb{R}^m$,
\begin{equation}
\left\{ {\omega \in \Omega :X\left(\omega\right) \bigcap K \ne \phi
} \right\} \in \mathfrak{F} \ .
\end{equation}
\end{definition}
In the sequel we assume that $X$ is discretely distributed, namely,
$X(\omega) \in \left\{A_0,...,A_n\right\}$, with $\Pr\left\{X = A_i
\right\}= \alpha_i \geq 0$ and $\sum\limits_{i = 0}^n {\alpha _i } =
1$. Moreover, we assume that $A_i \in \mathfrak{J}$, $i = 0,...n$.
Note that for any $f : \mathcal{F} \to \mathbb{R}$ and any random
set $X$, $f\left( X \right)$ defines a real-valued random variable.

Random closed sets $X_1,...,X_n$ are said to be \emph{independent}
if,
\begin{equation}\label{Independence}
\Pr \left\{ X_1 \in {\mathcal{X}}_{1},...,X_n  \in {\mathcal{X}}_{n}
\right\} = \Pr \left\{ {X_1  \in {\mathcal{X}}_{1} } \right\} \cdots
\Pr \left\{ {X_n \in {\mathcal{X}}_{n} } \right\} \ ,
\end{equation}
for all ${\mathcal{X}}_1,...,{\mathcal{X}}_n \in
\mathfrak{B}\left(\mathcal{F}\right)$. Here
$\mathfrak{B}\left(\mathcal{F}\right)$ is generated by all
collections of closed sets of the form $ \left\{ F \in {\mathcal{F}}
: F \bigcap K \ne \phi \right\}$ with $K$ running through all
compact subsets of $\mathbb{R}^m$ (\cite{molchanov2005theory},
Section 1.2).

The \emph{coverage function} $p_X \left( \cdot \right) :
\mathbb{R}^m \to [0,1]$ of the closed random set $X$ is
(\cite{molchanov2005theory}, Section 2.2)
\begin{equation}\label{CoverageFunction}
p_X \left( u \right) = \Pr \left\{ {u \in X} \right\} \ .
\end{equation}
Notice that for discretely distributed random set $X$,
\begin{equation}\label{DescreteCoverageFunction}
p_X \left( u \right) = \sum\limits_{\left\{ {i:u \in A_i } \right\}}
{\alpha _i } \ .
\end{equation}
The following relation is useful,
\begin{equation}\label{Robbin}
\int\limits_{\mathbb{R}^m } {p_X \left( u \right)du = E\left(\mu
\left( X \right)\right)} \ ,
\end{equation}
where $E$ denotes the \emph{expectation} of a real-valued random
variable. Clearly, for a discretely distributed $X$ accepting values
$\left\{A_0,...,A_n\right\}$ the integral in (\ref{Robbin}) can be
taken over $\bigcup\limits_{i = 0}^n {A_i } $.

Let $u \in \mathbb{R}^m$, setting $ {\mathcal{X}}_i  = \left\{ {F
\in {\mathcal{F}}: F \bigcap \left\{ u \right\} \ne \phi }
\right\}$, $i = 0,...,n$ one obtains from (\ref{Independence}) that
for independent $X_1,...X_n$,
\begin{equation}
\Pr \left\{ {u \in X_1 ,...,u \in X_n } \right\} = p_{X_1 } \left( u
\right)\cdots p_{X_n } \left( u \right) \ .
\end{equation}

In the next section we define a new average of sets, with which we
adapt the Bernstein operators to SVFs and obtain results analogous
to Theorems \ref{Weierstrass} and \ref{Kac}.

\section{The partition average of
sets}\label{sectionPartitionAverage} The construction of our average
of sets is built upon several definitions. We begin with
\begin{definition}\label{SubSetGenFunction}
Let $\Psi : \mathbb{J} \times \left[ {0,1} \right]  \to
\mathfrak{J}$ be such that
\begin{enumerate}
\item $\Psi \left( {A,t} \right) \subseteq A$
\item $ \label{SubSetMeasureProp} \mu \left( {\Psi \left( {A,t} \right)} \right) = t\mu \left( A \right)$
\item For $s \leq t$, $\Psi \left( {A,s} \right) \subseteq \Psi \left( {A,t} \right)$ .
\end{enumerate}
The function $\Psi $ is called the \emph{subset-generating
function}.
\end{definition}
Note that since $A \in \mathbb{J}$ and $\Psi \left( {A,t} \right)
\in \mathfrak{J}$,
\begin{equation}\label{ciOne}
\Psi \left( {A,0} \right) = \phi \ , \ \Psi \left( {A,1} \right) =
\operatorname{ci}\left(A \right) \ .
\end{equation}

\begin{figure}
\begin{center}
\includegraphics[scale=0.6]{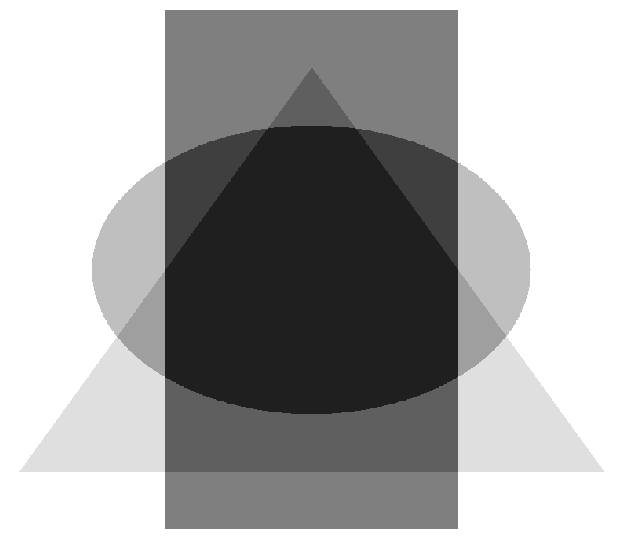}
\end{center}
\caption{The partition of the union of three subsets of
$\mathbb{R}^2$, the triangle, the rectangle and the ellipse. Regions
with similar gray-tone belong to the same element of the partition
of the union.} \label{fig:Partition}
\end{figure}

For any collection of sets in $\mathfrak{J}$, we consider a special
partition of their union to mutually disjoint sets,
\begin{definition}\label{PartionOfUnion} Let $ \left\{ {A_0 ,...,A_n } \right\} \subset \mathfrak{J}$.
For any subset $\chi$ of the indices $\left\{0,...,n\right\}$,  we
define the set
\begin{equation}\label{OmegaChiEq}
\Omega _\chi ^{A_0 ,...,A_n }  = \left( {\bigcap\limits_{k \in \chi
} {A_k } } \right)\backslash \left( {\bigcup\limits_{l \in \{
0,..,n\} \backslash \chi } {A_l } } \right) \ .
\end{equation}
For a fixed collection of sets $\left\{A_0,...,A_n\right\}$ we use
the shorthand notation $\Omega_\chi$. The collection of sets
\begin{equation*}
\left\{\Omega _\chi ^{A_0 ,...,A_n }: \chi \in 2^{\{n\}} \right\} \
,
\end{equation*}
where $2^{\{n\}}$ denotes all subsets of the set of integers
$\{0,...,n\}$, is termed the \emph{partition of the union} of
$A_0,...,A_n$. The sets  $\Omega_\chi$ are termed \emph{elements of
the partition}.
\end{definition}
An example of the partition of the union of three subsets on
$\mathbb{R}^2$ is given in Figure \ref{fig:Partition}.

Here we state several properties of the partition of the union, that
follow easily from Definition \ref{PartionOfUnion}.
\begin{lemma}\label{PropertiesOfThePartion} Let $ \left\{ {A_0 ,...,A_n } \right\}
\subset \mathfrak{J}$ then
\begin{enumerate}
\item \label{JordanElements} $\Omega _{\chi } \in \mathbb{J}$

\item \label{PartitionDisjoint}$\Omega _{\chi _{_1 } }  \bigcap \Omega _{\chi _2 }= \phi$, $\chi_1, \chi_2 \in 2^{\{n\}}$, $\chi_1 \ne \chi_2$.

\item\label{PartitionConstitute}$\bigcup\limits_{\left\{ {\chi  \in 2^{\{ n\} } :j \in \chi } \right\}} {\Omega _\chi   = A_j } $ .

\item \label{FixedTilde} For a fixed $\tilde \chi  \in 2^{\left\{ n
\right\}}$, $\tilde \chi \neq \phi$, $\bigcup\limits_{\left\{ {\chi
\in 2^{\left\{ n \right\}} :\tilde \chi  \subseteq \chi } \right\}}
{\Omega _\chi  } = \bigcap\limits_{j \in \tilde \chi } {A_j }$.

\item\label{PatitionEqualsUnion} $\bigcup\limits_{\chi  \in 2^{\{ n\} } } {\Omega _\chi  = } \bigcup\limits_{i = 0}^n {A_i } $

\end{enumerate}
\end{lemma}
Next observation connects the notion of the partition of union with
random sets. The proof of this observation follows from Definition
\ref{PartionOfUnion} and (\ref{DescreteCoverageFunction}).
\begin{lemma}\label{PartitionCoverage} Let $X$ be a random set, $X\left( \omega  \right) \in \left\{ {A_0 ,...,A_n } \right\}$,
$\Pr \left( {X = A_i } \right) = \alpha _i $. The coverage function,
$p_X \left( u \right)$ is constant over each element $\Omega _\chi $
of the partition of the union of $A_0,...,A_n$, and
\begin{equation*}
p_X \left( u \right)_{|\Omega _\chi  }  = \sum\limits_{i \in \chi }
{\alpha _i } \ .
\end{equation*}
\end{lemma}

We are now in a position to define a new weighted average of sets in
$\mathfrak{J}$, which is based on the partition of the union of the
averaged sets.
\begin{definition}\label{partitionAvearge}
Let $A_0,...,A_n \in \mathfrak{J}$  and $\alpha_0,...,\alpha_n \in
[0,1]$, $ \sum\limits_{i = 0}^n {\alpha _i  = 1} $. The
\emph{partition average} of $A_0,...,A_n$ with the weights
$\alpha_1,...,\alpha_n$ is
\begin{equation}\label{partitionAveargeEq}
\mathop  \otimes \limits_{i = 0}^n \alpha _i A_i :
 = \bigcup\limits_{\chi  \in 2^{\{ n\} } } {\Psi \left( {\Omega _\chi  ,\sum\limits_{k \in \chi } {\alpha _k } } \right)}
 \ ,
\end{equation}
where $\Psi$ is a subset-generating function in Definition
\ref{SubSetGenFunction}.
\end{definition}

Using the partition average, we can define expectation of a
discretely distributed random set as
\begin{definition}\label{partitionExpectation}
Let $X$ be a random set, $X\left( \omega  \right) \in \left\{ {A_0
,...,A_n } \right\}$,  $ \Pr \left( {X = A_i } \right)= \alpha _i ,i
= 0,...n$. The \emph{partition expectation} of $X$ is
\begin{equation}\label{partitionExpectationEq}
E_P\left(X\right): = \mathop  \otimes \limits_{i = 0}^n \alpha _i
A_i \ .
\end{equation}
\end{definition}
\begin{remark} In view of Lemma \ref{PartitionCoverage}, the
partition expectation is related to the coverage function through
\begin{equation}
E_P \left( X \right) = \bigcup\limits_{\chi  \in 2^n } {\Psi \left(
{\Omega _\chi  ,\left. {p_X } \right|_{\Omega _\chi  } } \right)} \
.
\end{equation}
\end{remark}
Next we state relevant properties of the partition average of sets.
\begin{theorem}\label{partitionAveargeProps}
In the notation of Definition \ref{partitionAvearge},
\begin{enumerate}
\item\label{partAverageClosureProp} $\mathop  \otimes \limits_{i = 0}^n \alpha _i A_i \in \mathfrak{J}$

\item\label{partAverageCommutProp} For any permutation $r(\cdot)$ of $\{0,...,n\}$, $\mathop  \otimes \limits_{i = 0}^n \alpha _i A_i  = \mathop  \otimes
\limits_{i = 0}^n \alpha _{r\left( i \right)} A_{r\left( i\right)}$

\item\label{partAverageAllSameProp}If for some $k \in\left\{1,...,n\right\}$, $A_k  = A_{k + 1}  = ... = A_n$, then $\mathop  \otimes
\limits_{i = 0}^n \alpha _i A_i  = \mathop  \otimes \limits_{i =
0}^k \beta _i A_i $ with $\beta_i = \alpha_i, i =0,...,k-1$ and
$\beta _k = \sum\limits_{i = k}^n {\alpha _i } $. In particular,
$\mathop \otimes \limits_{i = 0}^n \alpha _i A = A$
\item\label{partAverageInclusionProp} $
\operatorname{ci} \left( {\bigcap\limits_{\left\{ {i:\alpha _i  > 0}
\right\}}^n {\alpha _i A_i } } \right) \subseteq \mathop  \otimes
\limits_{i = 0}^n \alpha _i A_i  \subseteq \bigcup\limits_{\left\{
{i:\alpha _i  > 0} \right\}}^n {A_i }$

\item\label{partAverageInterpolationProp} If for some $j$, $\alpha _j  = 1$, then $\mathop  \otimes \limits_{i = 0}^n \alpha _i A_i  = A_j $

\item\label{partAverageMeasureProp} $\mu \left( {\mathop  \otimes \limits_{i = 0}^n \alpha _i A_i }
\right) = \sum\limits_{i = 0}^n {\alpha _i } \mu \left( {A_i
}\right)$

\end{enumerate}
\end{theorem}
\proof To obtain Property \ref{partAverageClosureProp}, observe that
by Definition \ref{SubSetGenFunction},  $\Psi \left( {\Omega _\chi
,t} \right) \in \mathfrak{J}$ for any $\chi  \in 2^{\{n\} }$, $t \in
[0,1]$, and recall that $\mathfrak{J}$ is closed under finite
unions. Properties
\ref{partAverageCommutProp},\ref{partAverageAllSameProp},
 follow immediately from the definition
of the partition average.

Next we prove Property \ref{partAverageInclusionProp}. Let $\tilde
\chi  = \left\{ {j \in \left\{ {0,...,n} \right\}:\alpha_j  > 0}
\right\}$. Since $\sum\limits_{j \in \tilde \chi } {\alpha _j  =
1}$, $\sum\limits_{j \in \chi } {\alpha _j  = 1}$ for $\chi
\supseteq \tilde \chi$. Therefore, from (\ref{ciOne}),
(\ref{ciOfUnion}) and Property \ref{FixedTilde} in Lemma
\ref{PropertiesOfThePartion},
\begin{equation*}
\bigcup\limits_{\left\{ {\chi :\tilde \chi  \subseteq \chi }
\right\}} {\Psi \left( {\Omega _\chi  ,\sum\limits_{i \in \chi }
{\alpha _i } } \right)}  = \bigcup\limits_{\left\{ {\chi :\tilde
\chi  \subseteq \chi } \right\}} {{\text{ci}}\left( {\Omega _\chi  }
\right)}  = {\text{ci}}\left( {\bigcup\limits_{\left\{ {\chi :\tilde
\chi  \subseteq \chi } \right\}} {\Omega _\chi  } } \right) =
\operatorname{ci} \left( {\bigcap\limits_{i \in \tilde \chi } {A_i }
} \right) \ ,
\end{equation*}
and thus $ \operatorname{ci} \left( {\bigcap\limits_{i \in \tilde
\chi } {A_i } } \right) \subseteq \mathop  \otimes \limits_{i = 0}^n
\alpha _i A_i$. The other part of Property
\ref{partAverageInclusionProp}, follows from the observation that $
\mathop  \oplus \limits_{i = 0}^n \alpha _i A_i  \subseteq
\bigcup\limits_{\left\{ {\chi :\chi  \cap \tilde \chi  \ne \phi }
\right\}} {\Omega _\chi  }  \subseteq \bigcup\limits_{j \in \tilde
\chi } {A_j }$.

Property \ref{partAverageInterpolationProp} is an immediate
consequence of Property \ref{partAverageInclusionProp}. Next we
prove Property \ref{partAverageMeasureProp}. From the definition of
the partition average, from the fact that the sets $ \left\{ {\Omega
_\chi  :\chi  \in 2^{\left\{ n \right\}} } \right\}$ are pairwise
disjoint and from the properties of the subset-generating function,
we obtain that
\begin{equation}\label{expectMes1}
\mu \left( {\mathop  \otimes \limits_{n = 1}^n \alpha _i A_i
}\right) = \mu \left( {\bigcup\limits_{\chi  \in 2^{\left\{ n
\right\}} } {\Psi \left( {\Omega _\chi  ,\sum\limits_{i \in \chi }
{\alpha _i } } \right)} } \right) = \sum\limits_{\chi  \in
2^{\left\{ n \right\}} } {\left( {\sum\limits_{i \in \chi } {\alpha
_i } } \right)\mu \left( {\Omega _\chi  } \right)} \ .
\end{equation}
To proceed with the proof of Property \ref{partAverageMeasureProp},
 we interpret $\mathop
\otimes \limits_{i = 0}^n \alpha _i A_i $ as the partition
expectation of a random set $X$, such that $\Pr \left( {X = A_i }
\right) = \alpha _i $, $i = 0,...,n$. Now by Lemma
\ref{PartitionCoverage} and by Properties \ref{PartitionDisjoint},
\ref{PatitionEqualsUnion} of Lemma \ref{PropertiesOfThePartion},
\begin{eqnarray}\label{expectMes2}
\sum\limits_{\chi  \in 2^{\left\{ n \right\}} } {\left(
{\sum\limits_{i \in \chi } {\alpha _i } } \right)\mu \left( {\Omega
_\chi  } \right)}  = \sum\limits_{\chi  \in 2^{\left\{ n \right\}} }
{p_X \left| {_{\Omega _\chi  } \mu \left( {\Omega _\chi  } \right)}
\right.}  = \\ \nonumber = \sum\limits_{\chi  \in 2^{\left\{ n
\right\}} } {\int\limits_{\Omega _\chi  } {p_X \left( u \right)du} }
= \int\limits_{\bigcup\limits_{i = 0}^n {A_i } } {p_X \left( u
\right)du}  \ .
\end{eqnarray}
Finally, we  apply (\ref{Robbin}) to obtain that
\begin{equation*}
\mu \left( {\mathop  \otimes \limits_{n = 1}^n \alpha _i A_i }
\right) = E\left(\mu \left( X \right)\right) = \sum\limits_{i = 0}^n
\alpha_i {\mu \left( {A_i } \right)} \ .
\end{equation*}
\qed

\begin{remark}\label{Analogy} The above properties of the partition average are analogous to those of weighted averages between non-negative numbers.
In this analogy, the measure of a set replaces the absolute value of
a number, the measure of the symmetric difference of two sets
$(d_\mu (\cdot; \cdot))$ replaces the absolute value of the
difference between two numbers. Moreover, the intersection and union
of sets replace the minimum and the maximum of numbers, and finally
the relation $\subseteq$ between sets replaces the relation $\leq$
between numbers.
\end{remark}

The next theorem treats the distance between the partition
expectations of two independent random sets distributed over the
same collection of sets $\left\{ {A_0 ,...,A_n } \right\}$.
\begin{theorem}\label{ExpectOfDistances}
Let $X_1$,$X_2$ be independent random sets, $\Pr \left\{ {X_1  = A_i
} \right\} = \alpha _i $, $\Pr \left\{ {X_2  = A_i } \right\} =
\beta _i $, $i = 0,...,n$, with $\sum\limits_{i = 0}^n {\alpha _i  =
\sum\limits_{i = 0}^n {\beta _i  = 1} } $.
 Then
\begin{equation*}
d_\mu  \left( {E_P \left(X_1\right) ,E_P \left(X_2\right) } \right)
\leq E\left(d_\mu \left( {X_1 ,X_2 } \right)\right) \ ,
\end{equation*}
where $d_\mu \left( {X_1 ,X_2 } \right)$ is the real-valued random
variable $d_\mu  \left( {X_1 ,X_2 } \right) = \mu \left( {X_1 \Delta
X_2 } \right)$, namely
\begin{equation*}
\Pr \left( {d{}_\mu \left( {X_1 ,X_2 } \right) = \mu \left( {A_i
\Delta A_j } \right)} \right) = \alpha _i \beta j, \ i,j = 0,...,n \
.
\end{equation*}

\end{theorem}
\proof

It follows from the definition of the partition expectation, and by
the fact that the partition elements are disjoint sets with their
union equal to $ {\bigcup\limits_{i = 0}^n {A_i } }$ (see Lemma
\ref{PropertiesOfThePartion}), that
\begin{equation*}
d_\mu  \left( {E_P \left( {X_1 } \right),E_P \left( {X_2 } \right)}
\right) = \sum\limits_{\chi  \in 2^{\left\{ n \right\}} } {d_\mu
\left( {\Psi \left( {\Omega _\chi ^{} ,\sum\limits_{i \in \chi }
{\alpha _i } } \right),\Psi \left( {\Omega _\chi ,\sum\limits_{i \in
\chi } {\beta _i } } \right)} \right)}   \ .
\end{equation*}
By the properties of the subset-generating function, for any $\chi
\in 2^{\left\{ n \right\}} $ one of the two sets $\Psi \left(
{\Omega _\chi  ,\sum\limits_{i \in \chi } {\alpha _i } }
\right),\Psi \left( {\Omega _\chi ,\sum\limits_{i \in \chi } {\beta
_i } } \right)$ is necessarily contained in the other, and we get
from (\ref{IncludedDistance}),
\begin{equation*} d_\mu  \left( {E_P \left(X_1\right)
,E_P \left(X_2\right) } \right) = \sum\limits_{\chi  \in 2^{\left\{
n \right\}} } {\mu \left( {\Omega _\chi} \right)\left|
{\sum\limits_{i \in \chi } {\alpha _i }  - \sum\limits_{i \in \chi }
{\beta _i } } \right|} \ .
\end{equation*}
Now from Lemma \ref{PartitionCoverage} and Properties
\ref{PartitionDisjoint},\ref{PatitionEqualsUnion} of Lemma
\ref{PropertiesOfThePartion} we get
\begin{equation}\label{DistanceOfExpects}
d_\mu  \left( {E_P \left(X_1\right) ,E_P \left(X_2\right) } \right)
= \int\limits_{\bigcup {A_i } } {\left| {p_{X_1 } \left( u \right) -
p_{X_2 } \left( u \right)} \right|} du \ .
\end{equation}
On the other hand,
\begin{equation*}
E\left(d_\mu  \left( {X_1 ,X_2 } \right)\right) = E\left( {\mu
\left( {X_1 \Delta X_2 } \right)} \right) = E\left( {\mu \left( {X_1
\backslash X_2 } \right)} \right) + E\left( {\mu \left( {X_2
\backslash X_1 } \right)} \right) \ .
\end{equation*}
Since $X_1, X_2$ are independent,
\begin{equation*}
p_{_{X_1 \backslash X_2 } } \left( u \right) = p_{X_1 } \left( u
\right)\left( {1 - p_{X_2 }\left(u\right) } \right) , \ u \in
\mathbb{R}^m \ ,
\end{equation*}
and we obtain from (\ref{Robbin}),
\begin{equation*}
E\left( {\mu \left( {X_1 \backslash X_2 } \right)} \right) =
\int\limits_{\bigcup\limits_{i = 0}^n {A_i } } {p_{X_1 } \left( u
\right)\left( {1 - p_{X_2 } \left( u \right)} \right)du}  \ .
\end{equation*}
Using similar observations for $ E\left( {\mu \left( {X_2 \backslash
X_1 } \right)} \right)$, we arrive at
\begin{equation}\label{ExpectOfDistances}
E\left( {d_\mu  \left( {X_1 ,X_2 } \right)} \right) =
\int\limits_{\bigcup\limits_{i = 0}^n {A_i } } {\left[ {p_{X_1 }
\left( u \right)\left( {1 - p_{X_2 } \left( u \right)} \right) +
p_{X_2 } \left( u \right)\left( {1 - p_{X_1 } \left( u \right)}
\right)} \right]du}  \ .
\end{equation}
It is easy to obtain the claim of the theorem, by inspecting the
relations (\ref{DistanceOfExpects}) and (\ref{ExpectOfDistances}),
since
\begin{equation}\label{algebraicRelation}
\left| {a - b} \right| \leqslant a\left( {1 - b} \right) + b\left(
{1 - a} \right), \ a,b \in [0,1] \ .
\end{equation}
 \qed

From the above theorem we obtain,
\begin{corollary}\label{ExpectationToOne} Let $X$ be a random set, $\Pr \left\{ {X  = A_i } \right\} = \alpha
_i$, $i = 0,...,n$. Then
\begin{equation}\label{ExpectOfDistancesToOne}
d_\mu  \left( {E_P\left( X \right),A_j } \right) = E\left(d_\mu
\left( {X,A_j } \right)\right) \ ,
\end{equation}
for any $j \in \left\{0,...,n\right\}$.
\end{corollary}
\proof Consider the random set $\widetilde X$, $\Pr \left\{
{\widetilde X = A_i } \right\} = \delta _{ij}$, $i = 0,..,n$ with
$\delta _{ij} = 1$ for $i = j$ and $\delta _{ij} = 0$ otherwise. By
Property \ref{partAverageInterpolationProp} of Theorem
\ref{partitionAveargeProps}, $E_P \left(\widetilde X \right) = A_j$,
thus from Theorem \ref{ExpectOfDistances},
\begin{equation}\label{ExpectToOneEx}
 d_\mu  \left( {E_P
\left(X\right) , A_j} \right) \leq Ed_\mu \left( {X ,A_j } \right) \
.
\end{equation}
Since $ p_{\widetilde X } \left( u \right) \in \left\{ {0,1}
\right\}$, if follows from (\ref{algebraicRelation}), that there is
an equality in (\ref{ExpectToOneEx}). \qed

Approximation results in the next section are based upon the
following corollary, which is derived from Corollary
\ref{ExpectationToOne}.
\begin{corollary}\label{AverageToOne}Let $A_0,...,A_n \in \mathfrak{J}$  and $\alpha_0,...,\alpha_n \in
[0,1]$, $ \sum\limits_{i = 0}^n {\alpha _i  = 1} $, then
\begin{equation}\label{AverageOfDistancesToOne}
d_\mu  \left( {A_j ,\mathop  \otimes \limits_{i = 0}^n \alpha _i A_i
} \right) = \sum\limits_{i = 0}^n {\alpha _i d_\mu  \left( {A_j ,A_i
} \right)} \ .
\end{equation}
Namely, the distance to the partition average from any of the
averaged sets is equal to the average of the distances from this set
to all the averaged sets.
\end{corollary}

\begin{remark}\label{analogyDistance}
Relation (\ref{AverageOfDistancesToOne}) is reminiscent of the
relation,
\begin{equation}\label{ExpectOfDistancesToOneNumbers}
\left| {r - \sum\limits_{i = 0}^n {\alpha _i p_i } } \right|
\leqslant \sum\limits_{i = 0}^n {\alpha _i \left| {r - p_i } \ ,
\right|}  \ .
\end{equation}
where $r,p_i$, $i=0,...,n$ are numbers. However, notice that there
is equality in (\ref{AverageOfDistancesToOne}) versus inequality in
(\ref{ExpectOfDistancesToOneNumbers}). Moreover, observe that
(\ref{AverageOfDistancesToOne}) applies only to the sets
participating in the partition average, while
(\ref{ExpectOfDistancesToOneNumbers}) applies to any $r \in
\mathbb{R}$. This limitation has implications to the approximation
power of methods based on the partition average.
\end{remark}

The partition average of two sets possesses also the \emph{metric
property} \cite{dyn2001spline} relative to
$d_\mu\left(\cdot,\cdot\right)$.
\begin{corollary}\label{metricProperty}
Let $A_0,A_1 \in \mathfrak{J}$, $\alpha_0,\beta_0 \in [0,1]$, then
\begin{equation}
d_\mu  \left( {\alpha _0 A_0  \otimes \left( {1 - \alpha _0 }
\right)A_1 ,\beta _0 A_0  \otimes \left( {1 - \beta _0 } \right)A_1
} \right) = \left| {\alpha _0  - \beta _0 } \right|d_\mu  \left(
{A_0 ,A_1 } \right) \ .
\end{equation}
\end{corollary}
\proof Let $X_1,X_2$ be random sets such that $ \Pr \left( {X_1  =
A_i } \right) = \alpha _i$ and $ \Pr \left( {X_2  = A_i } \right) =
\beta _i$, $i = 0,1$. Note that $\alpha _1  = 1 - \alpha _0 ,\beta
_1  = 1 - \beta _0$. By (\ref{DistanceOfExpects}),
\begin{equation*}
d_\mu  \left( {E_P \left( {X_1 } \right),E_P \left( {X_2 } \right)}
\right) = \int\limits_{A_0  \cup A_1 } {\left| {p_{X_1 } \left( u
\right) - p_{X_2 } \left( u \right)} \right|} du
 \ .
\end{equation*}
Since for $u \notin A\Delta B$, $ p_{X_1 } \left( u \right) = p_{X_2
} \left( u \right)$, we get
\begin{equation*}
d_\mu  \left( {\alpha _0 A_0  \otimes \left( {1 - \alpha _0 }
\right)A_1 ,\beta _0 A_0  \otimes \left( {1 - \beta _0 } \right)A_1
} \right) = \int\limits_{A\Delta B} {\left| {\alpha _0  - \beta _0 }
\right|} du = \left| {\alpha _0  - \beta _0 } \right|\mu \left(
{A\Delta B} \right) \ .
\end{equation*}
\qed

\begin{figure}
\begin{center}
\includegraphics[scale=0.6]{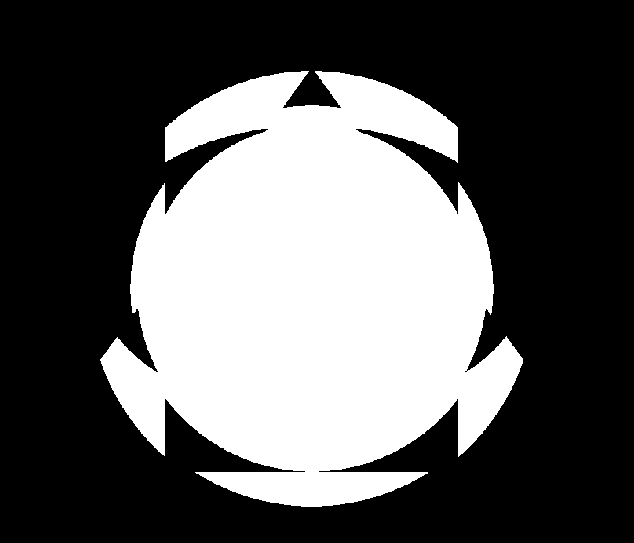}
\end{center}
\caption{The partition average $\mathop  \otimes \limits_{i = 0}^2
\alpha _i A_i $, with $\alpha _i = \frac{1}{3}, i = 0,1,2$ and $A_i,
i = 0,1,2$ the three sets in Figure \ref{fig:Partition}. }
\label{fig:AverageExample}
\end{figure}

To complete the construction of the partition average of sets, we
need to provide a concrete example of a subset-generating function
in Definition \ref{SubSetGenFunction}. We denote by $Bl\left( {p,r}
\right) $ a ball of radius $r$ about $p \in \mathbb{R}^m$, namely
\begin{equation*}
Bl\left( {p,r} \right) = \left\{ {q \in \mathbb{R}^m :\left\| {q -
p} \right\| \leqslant r} \right\} \ ,
\end{equation*}
with $\|\cdot\|$ the Euclidean norm on $\mathbb{R}^m$. The
subset-generating function $\Psi :\mathbb{J} \times [0,1] \to
\mathfrak{J}$ is defined by
\begin{equation}\label{Realization}
\Psi \left( {A,t} \right) = \operatorname{ci} \left( {Bl\left(
{p,r_{A,t} } \right)\bigcap A } \right) \ ,
\end{equation}
where $r_{A,t}$ is chosen so that $ \mu \left( {\Psi \left( {A,t}
\right)} \right) = t\mu \left( A \right)$. The existence of $r_{A
,t}$ as above for any $t \in [0,1]$ follows from the continuity of
the volume of the ball as a function of its radius. An example of
the partition average with a such defined subset-generating function
$\Psi$ is shown in Figure \ref{fig:AverageExample}. In this example,
$p$ is the \emph{centroid} of the union of the averaged sets.

Although there is a significant resemblance between the partition
average of sets and the weighted average of numbers as is noticed in
Remarks \ref{Analogy} and \ref{analogyDistance}, the partition
average of sets lacks several important properties of the weighted
average of numbers.
\begin{remark}\label{remarkNotAssoc}The partition average is generally not associative,
\begin{equation*}
\begin{array}{*{20}c}
 \mathop  \otimes \limits_{i = 0}^2 \alpha _i A_i &  \ne &
\left( {\alpha _0  + \alpha _1 } \right)\left( {\frac{{\alpha _0 }}
{{\alpha _0  + \alpha _1 }}A_0  \otimes \frac{{\alpha _1 }} {{\alpha
_0  + \alpha _1 }}A_1 } \right) \otimes \alpha _2 A_2  \\
& \ne & a_0 A_0 \otimes \left( {\alpha _1  + \alpha _2 }
\right)\left( {\frac{{\alpha _1 }} {{\alpha _1  + \alpha _2 }}A_1
\otimes \frac{{\alpha _2 }} {{\alpha _1  + \alpha _2 }}A_2 } \right)
\ .
\end{array}
\end{equation*}
\end{remark}
\begin{remark}
Zero-weighted sets in (\ref{partitionAveargeEq}) affect the
partition average by affecting the partition of the union of all the
sets, namely,
\begin{equation*}\label{remarkShort}
\mathop  \otimes \limits_{i = 0}^n \alpha _i A_i  \ne \mathop
\otimes \limits_{i = 0}^{n + 1} \alpha _i A_i \ ,
\end{equation*}
with $\sum\limits_{i = 0}^n {\alpha _i }  = 1$ and $\alpha _{n + 1}
= 0$. Yet, the average "is between the intersection and the union of
the sets with positive weights" (see Property
\ref{partAverageInclusionProp} in Theorem
\ref{partitionAveargeProps}).
\end{remark}

\section{Set-valued Bernstein approximation based on the partition average}\label{sectionBernstein}

Using the partition average of sets defined in the previous section,
we can now define a set-valued operator analogous to
(\ref{realBernstein}),
\begin{definition}\label{setBernstein}
The the set-valued Bernstein operator is the mapping
$F\left(\cdot\right) \to B_n\left(F, \cdot\right)$ given by
\begin{equation}\label{setBernstein}
B_n \left( {F,x} \right) = \mathop  \otimes \limits_{i = 0}^n
b\left( {n,x;i} \right)F\left( {\frac{i} {n}} \right), \; x \in
[0,1]  \; ,
\end{equation}
for any $F:\left[ {0,1} \right] \to \mathfrak{J}$, where
$b\left({n,x;i} \right)$ are defined in (\ref{Bn}).
\end{definition}

Using Definition \ref{setBernstein}, we aim to obtain approximation
results analogous to Theorem \ref{Weierstrass} and Theorem
\ref{Kac}. First we note that by Property
\ref{partAverageInterpolationProp} in Theorem
\ref{partitionAveargeProps} and by (\ref{Bn}), $ B_n \left( {F,0}
\right) = F\left( 0 \right),B_n \left( {F,1} \right) = F\left( 1
\right) $. The set-valued version of Theorem \ref{Weierstrass} is
\begin{theorem}\label{ConinuousApproximation}
Let $F : [0,1] \to \mathfrak{J}$ be a continuous SVF, then for any
$\varepsilon > 0$ there exists $N > 0$, such that for all $n \geq N$
and all $x \in [0,1]$,
\begin{equation}\label{ConinuousApproximationEq}
d_\mu  \left( {F\left( x \right),B_n \left( {F,x} \right)} \right) <
\varepsilon \ .
\end{equation}
\end{theorem}

The proof of Theorem \ref{ConinuousApproximation} is based on the
following two lemmas.
\begin{lemma}\label{LemmaClosetPoint}
In the notation of Definition \ref{setBernstein}, let $x \in [0,1]$
and let $x^{\prime}$ be the point closest to $x$ among
$\frac{i}{n}$, $i = 0,...,n$. Then
\begin{equation}
d_\mu  \left( {F\left( x \right),B_n \left( {F,x} \right)} \right)
\leq 2d_\mu  \left( {F\left( {x^\prime  } \right),F\left( x \right)}
\right) + \sum\limits_{i = 0}^n {b\left( {n,x;i} \right)d_\mu \left(
{F\left( x \right),F\left( {x_i } \right)} \right)}  \ .
\end{equation}
\end{lemma}
\proof By the triangle inequality,
\begin{equation}\label{triang1}
d_\mu  \left( {F\left( x \right),B_n \left( {F,x} \right)} \right)
\leq d_\mu  \left( {F\left( x \right),F\left( {x^{\prime}} \right)}
\right) + d_\mu  \left( {F\left( {x^{\prime}} \right),B_n \left(
{F,x} \right)} \right)  \ .
\end{equation}
We obtain from Corollary \ref{AverageToOne} that
\begin{equation}\label{firstEquality}
d_\mu  \left( {F\left( {x^\prime  } \right),B_n \left( {F,x}
\right)} \right) = \sum\limits_{i = 0}^n {b\left( {n,x;i}
\right)d_\mu  \left( {F\left( {x'} \right),F\left( {x_i } \right)}
\right)} \ ,
\end{equation}
and by the triangle inequality,
\begin{equation*}
\sum\limits_{i = 0}^n {b\left( {n,x;i} \right)d_\mu  \left( {F\left(
{x'} \right),F\left( {x_i } \right)} \right)} \leq d_\mu  \left(
{F\left( {x^\prime  } \right),F\left( x \right)} \right) +
\sum\limits_{i = 0}^n {b\left( {n,x;i} \right)d_\mu  \left( {F\left(
x \right),F\left( {x_i } \right)} \right)} \ .
\end{equation*}
This together with (\ref{triang1}) and (\ref{firstEquality})
completes the proof. \qed

\begin{lemma}\label{ContinuousDistanceExpectation}
Let $F:\left[ {0,1} \right] \to \mathfrak{J}$ be a continuous
function, then for any $\varepsilon > 0$ there exists $N>0$, such
that for all $n \geq N$ and all $x \in [0,1]$ ,
\begin{equation}\label{ContinuousBoundEq}
 \sum\limits_{i = 0}^n {b\left(
{n,x;i} \right)d_\mu  \left( {F\left( x \right),F\left( {x_i }
\right)} \right)} < \varepsilon  \ .
\end{equation}
\end{lemma}
\proof For any $x \in [0,1]$, consider the function
$g_x\left(\cdot\right) : [0,1] \to \mathbb{R}$,
\begin{equation}\label{gxFunction}
g_x \left( y \right) = d_\mu  \left( {F\left( x \right),F\left( y
\right)} \right) \ .
\end{equation}
It is easy to observe that due to the continuity of $F$, the family
of functions $\left\{ {g_x \left(  \cdot  \right):x \in [0,1]}
\right\}$ is uniformly equicontinuous and uniformly bounded.

Next we note that $ \sum\limits_{i = 0}^n {b\left( {n,x;i}
\right)d_\mu  \left( {F\left( x \right),F\left( {x_i } \right)}
\right)}$ is $B_n \left( {g_x ,x} \right)$. Since $g_x\left(x\right)
= 0$, we use Theorem \ref{Weierstrass} to conclude that there exists
$N_x > 0$, such that (\ref{ContinuousBoundEq}) holds for all $n \geq
N_x$. By the uniform continuity and boundedness of $\left\{ {g_x }
\right\}_{x \in \left[ {0,1} \right]} $ there exists $N < \infty$
such that $N_x \leq N$ for all $x \in [0,1]$. \qed

\proof \emph{of Theorem \ref{ConinuousApproximation}} $\ $ By Lemma
\ref{ContinuousDistanceExpectation}, there exists a positive integer
$N_1$ s.t. for $n > N_1$,
\begin{equation*}
\sum\limits_{i = 0}^n {b\left( {n,x;i} \right)d_\mu  \left( {F\left(
x \right),F\left( {x_i } \right)} \right)} < \frac{\varepsilon}{2}
\;, \; x \in [0,1] \; .
\end{equation*}
By continuity of $F$, there exists a positive integer $N_2$ s.t. for
$\left| {x - y} \right| < \frac{1} {{N_2 }}$, $x,y \in [0,1]$,
$d_\mu \left( {F\left( x \right),F\left( y \right)} \right) <
\frac{\varepsilon } {4}$. We set $N = \max \left\{ {N_1 ,N_2 }
\right\}$ and apply Lemma \ref{LemmaClosetPoint} to obtain the claim
of the theorem. \qed

Next we consider the rate of approximation of H{\"o}lder continuous
SVFs by the set-valued Bernstein operator. The class of H{\"o}lder
continuous SVFs, $\operatorname{Lip}\left(L,\nu\right)$, is defined
as in (\ref{HolderFunctionEq}), using $d_\mu \left( { \cdot , \cdot
} \right)$ instead of the distance between numbers.

\begin{theorem}\label{HolderApproximation}
Let $F \in \operatorname{Lip}\left(L,\nu\right)$, then for all $n
\in \mathbb{N}$ and all $x \in [0,1]$,
\begin{equation}\label{setValuedKacEx}
d_\mu  \left( {F\left( x \right),B_n \left( {F,x} \right)} \right)
\leq L\left( {\frac{1} {n}} \right)^\nu   + L\left( {\frac{{x\left(
{1 - x} \right)}} {n}} \right)^{\nu /2} \ .
\end{equation}
\end{theorem}
\proof For any $x \in [0,1]$, let the function
$g_x\left(\cdot\right)$ be given by (\ref{gxFunction}). $F \in
\operatorname{Lip}\left(L,\nu\right)$ implies that $g_x \in
\operatorname{Lip}\left(L,\nu\right)$ for all $x \in [0,1]$. Since
$g_x\left(x\right) = 0$, we get from Theorem \ref{Kac},
\begin{equation}\label{LipBound}
\sum\limits_{i = 0}^n {b\left( {n,x;i} \right)d_\mu  \left( {F\left(
x \right),F\left( {x_i } \right)} \right)} = B_n \left( {g_x ,x}
\right) \le L\left( {\frac{{x\left( {1 - x} \right)}}{n}}
\right)^{\nu /2} \;.
\end{equation}
From $F \in \operatorname{Lip}\left(L,\nu\right)$, (\ref{LipBound})
and Lemma \ref{LemmaClosetPoint}, we obtain the result of the
theorem. \qed

Note that for a fixed $x \in \left(0,1\right)$, the term $L\left(
{\frac{1} {n}} \right)^\nu $ in (\ref{setValuedKacEx}) is dominated
by $L\left( {\frac{{x\left( {1 - x} \right)}} {n}} \right)^{\nu
/2}$, so the results obtained in Theorem \ref{Kac} and Theorem
\ref{HolderApproximation} are asymptotically equivalent.

\section{Set-valued Bernstein approximation with the de Casteljau's algorithm}\label{sectionCastelio}
A widely used method for the evaluation of the real-valued Bernstein
operators $B_n\left(f,x\right)$ is the de Casteljau's algorithm
(see, e.g. \cite{farin2002curves}, Chapter 4). The algorithm
evaluates $B_n\left(f,x\right)$ through a sequence of averages of
two numbers (\emph{binary averages}) and is based on the following
recurrence relation,
\begin{equation}\label{bernRealRecEq}
b\left( {n,x;i} \right) = \left( {1 - x} \right)b\left( {n - 1,x;i}
\right) + xb\left( {n - 1,x;i - 1} \right)
 \ ,
\end{equation}
where $b\left( {n,x;i} \right)$ are given in (\ref{Bn}). $B_n \left(
{f,x} \right)$ in (\ref{realBernstein}) can be represented using
(\ref{bernRealRecEq}) as,
\begin{equation}\label{bernRealBnk}
B_n \left( {f,x} \right) = \sum\limits_{i = 0}^n {b\left( {n,x;i}
\right)f_i^n  = \sum\limits_{i = 0}^{n - 1} {b\left( {n - 1,x;i}
\right)} } f_i^{n-1} \ ,
\end{equation}
with
\begin{equation}\label{bernRealFk}
f_i^n  = f\left( {\frac{i}{n}} \right), \ i = 0,...,n \ \text{and} \
f_i^{n-1} = \left( {1 - x} \right)f_i^{n}  + xf_{i + 1}^{n} \ , i =
0,1,...,n-1 \ .
\end{equation}
The de Casteljau's algorithm repeats this recursion $n$ times to get
\begin{equation}\label{BnRecFinal}
B_n(f,x) = f_0^0 \ .
\end{equation}
In the real-valued case, (\ref{realBernstein}) and the recursive
relations (\ref{bernRealBnk})-(\ref{BnRecFinal}) are equivalent,
though the evaluation of $B_n\left(f,x\right)$ by the de Casteljau's
algorithm is numerically stable.

A straightforward adaptation to SVFs of the recursive relations
(\ref{bernRealBnk})-(\ref{BnRecFinal}) based on the partition
average is
\begin{equation}
F_i^k  =  \left( {1 - x} \right)F_i^{k + 1}  \otimes xF_{i + 1}^{k +
1}, \ i = 0,...k, \ k = n-1,...,0\ ,
\end{equation}
with $F_i^n  = F\left( {\frac{i}{n}} \right)$, and $\hat B_n \left(
{F,x} \right)$ is set to be $F_0^0 $. Note that by Remark
\ref{remarkNotAssoc}, this adaptation yields a set-valued operator
which is different from that in (\ref{setBernstein}). The above
construction is similar to that in \cite{dyn2006approximations} with
the metric average.  Similarly to \cite{dyn2006approximations}, we
do not expect that for a general SVF F, ${\hat B_n \left( {F,x}
\right)}$ converges to $F\left(x\right)$ as $n \to \infty$.

We now alter the adaptation and apply Corollary \ref{AverageToOne},
to obtain approximation results similar to Theorems
\ref{ConinuousApproximation} and \ref{HolderApproximation} also in
the case of the de Casteljau's representation of the Bernstein
operators. To this end we use a binary average, based on the
partition determined by $F\left(\frac{i}{n}\right)$, $i = 0,...,n$,
of the form
\begin{equation}\label{tildeOtimes}
\lambda A\tilde  \otimes \left( {1 - \lambda } \right)B = \mathop
\otimes \limits_{i = 0}^{n + 2} \beta _i E_i \ ,
\end{equation}
where $E_i = F_i$, $i = 0,...,n$, $E_{n+1} = A$, $E_{n+2} = B$ and
$\beta _i =0$, $i = 0,...,n$, $\beta _{n+1} = \lambda$, $\beta
_{n+2} = 1 - \lambda$. Then we apply the de Casteljau's algorithm
with this average, namely
\begin{equation}\label{tildeReqursion}
F_i^k  =  \left( {1 - x} \right)F_i^{k + 1}  \tilde \otimes xF_{i +
1}^{k + 1}, \ i = 0,...k,  \ k = n-1,...,0, \,
\end{equation}
and set,
\begin{equation}\label{setBnRecusrive}
B^{DC}_n \left( {F,x}\right) = F_0^0 \ .
\end{equation}

Since ${F\left( {\frac{i}{n}} \right)}$, $i = 0,...,n$, is in the
partition behind $\tilde \otimes$, we get from
(\ref{tildeReqursion}) and Corollary \ref{AverageToOne}, that
\begin{equation*}
\begin{array}{*{20}l}
d_\mu  \left( {F_0^0 ,F\left( {\frac{i} {n}} \right)} \right) &=&
d_\mu  \left( {\left( {1 - x} \right)F_0^1 \tilde  \otimes xF_1^1
,F\left( {\frac{i} {n}} \right)} \right)\\
 &=& \left( {1 - x}
\right)d_\mu  \left( {F_0^1 ,F\left( {\frac{i} {n}} \right)} \right)
+ xd_\mu  \left( {F_1^1 ,F\left( {\frac{i} {n}} \right)} \right) \ .
\end{array}
\end{equation*}
Continuing the recursion we finally obtain the real-valued de
Casteljau's algorithm for the function $g_{\frac{i}{n}}$ defined in
(\ref{gxFunction}), namely with the initial data
\begin{equation*}
g_{\frac{i}{n}} \left( {\frac{j}{n}} \right) = d_\mu  \left(
{F\left( {\frac{i} {n}} \right),F\left( {\frac{j}{n}} \right)}
\right), \ j = 0,...n.
\end{equation*}
Therefore by (\ref{setBnRecusrive}) we have for $i = 0,...,n$,
\begin{equation}\label{deDistance}
d_\mu  \left( {B^{DC}_n \left( {F,x} \right),F\left( {\frac{i} {n}}
\right)} \right) = \sum\limits_{j = 0}^n {b\left( {n,x;j}
\right)d_\mu\left( {F\left( {\frac{i} {n}} \right),F\left( {\frac{j}
{n}} \right)} \right)} , x \in [0,1] \ .
\end{equation}
The equality (\ref{deDistance}) is the same as
(\ref{firstEquality}), and therefore Theorems
\ref{ConinuousApproximation} and \ref{HolderApproximation} also
apply to the set-valued Bernstein operators defined by the de
Casteljau's algorithm with the binary average (\ref{tildeOtimes}).

\section{Approximation of SVFs by positive sample-based operators}\label{sectionOperators}
We consider the adaptation to SVFs of families of positive
sample-based operators for real-valued functions, defined for $n \in
\mathbb{N}$ as
\begin{equation}\label{RealPositiveOperatorEq}
\widetilde O_n \left( {f,x} \right) = \sum\limits_{i = 0}^{l_n }
{c_{n,i} \left( x \right)f\left( {x_{n,i} } \right)}  \ ,
\end{equation}
where $f : [0,1]\to \mathbb{R}$, $x \in [0,1]$, $0 = x_{n,0} <
x_{n,1}<...<x_{n,l_{n}}=1$, $c_{n,i}(x) \geq 0$ and $\sum\limits_{i
= 0}^{l_n } {c_{n,i} \left( x \right) = 1} $. Moreover, we denote
$\delta_n \left(x\right) = \mathop {\min }\limits_{i \in \{
0,...,l_n \} } \left| {x - x_{n,i} } \right|$ and assume that for
any $x \in [0,1]$, $\mathop {\lim }\limits_{n \to \infty
}\delta_n\left(x\right) = 0$. The real-valued Bernstein
approximation operators are a prominent example of a family of
operators as above. Other examples are the piecewise linear
interpolation operator and the Schoenberg spline operators
\cite{rabut1992introduction}.

In analogy to the adaptation of the Bernstein operators to SVFs, we
define for $F : [0,1] \to \mathfrak{J}$,
\begin{equation}\label{PositiveOperatorEq}
O_n \left( {F,x} \right) = \mathop  \otimes \limits_{i = 0}^{l_n }
c_{n,i} \left( x \right)F\left( {x_{n,i} } \right)  \ .
\end{equation}

The first result is obtained immediately due to Property
\ref{partAverageMeasureProp} in Theorem \ref{partitionAveargeProps}
of the partition average,
\begin{corollary}\label{MeasureTransition}
Let $\widetilde O_n$  and $O_n$ be as in
(\ref{RealPositiveOperatorEq}) and (\ref{PositiveOperatorEq})
respectively. Then
\begin{equation*}
\mu \left( {O_n\left( {F,x} \right)} \right) = \widetilde O_n\left(
{\mu \left( {F\left( x \right)} \right),x} \right) \; , x \in [0,1]
\;.
\end{equation*}
\end{corollary}

Next we extend the results obtained in Theorems
\ref{ConinuousApproximation} and \ref{HolderApproximation} for the
Bernstein operators. Using the method of proof in Theorem
\ref{ConinuousApproximation} we obtain
\begin{corollary}\label{ContinuousApproxExtension}
Let $\widetilde O_n$  and $O_n$ be as in
(\ref{RealPositiveOperatorEq}) and (\ref{PositiveOperatorEq})
respectively. Assume that for any continuous real-valued function $f
: [0,1] \to \mathbb{R}$ and any $\epsilon > 0$, there exists
$\widetilde N_{f,\epsilon}$ such that for all $n \geq \widetilde
N_{f,\epsilon}$ and all $x \in [0,1]$,
\begin{equation*}
\left| {f\left( x \right) - \widetilde O_n\left( {f,x} \right)}
\right| < \epsilon \; .
\end{equation*}
Then for any continuous SVF, $F : [0,1] \to \mathfrak{J}$, and any
$\epsilon > 0$, there exists $N_{F,\epsilon}$ such that for all $n
\geq N_{F,\epsilon}$ and all $x \in [0,1]$,
\begin{equation*}
d_\mu  \left( {F\left( x \right),O_n\left( {F,x} \right)} \right) <
\epsilon \; .
\end{equation*}
\end{corollary}

For H{\"o}lder continuous SVFs, we obtain by arguments similar to
those in the proof of Theorem \ref{HolderApproximation},
\begin{corollary}\label{HolderApproximationEx} Let $\widetilde O_n$  and $O_n$ be as in
(\ref{RealPositiveOperatorEq}) and (\ref{PositiveOperatorEq})
respectively. Define the approximation error of $\widetilde O_n$ to
functions in $\operatorname{Lip}\left(L,\nu\right)$ at $x \in [0,1]$
as
\begin{equation*}
e_{n,L,\nu } \left( x \right) = \mathop {\sup }\limits_{f \in
{\mathop{\rm Lip}\nolimits} \left( {\nu ,L} \right)} \left|
{\widetilde O_n \left( {f,x} \right) - f\left( x \right)} \right| \;
.
\end{equation*}
Then,
\begin{equation*}
\mathop {\sup }\limits_{F \in {\mathop{\rm Lip}\nolimits} \left(
{\nu ,L} \right)} d_\mu  \left( {O_n \left( {F,x} \right),F\left( x
\right)} \right) \le e_{n,L,\nu } \left( x \right) + L\delta _n
\left( x \right)^\nu  \; .
\end{equation*}
\end{corollary}

Approximation of continuous functions by positive operators are
discussed in (\cite{fellerintroduction}, Chapter VII, \S 1) in the
context of probability theory, while approximation results for
H{\"o}lder continuous real-valued functions are the subject of
\cite{mathe2003asymptotic}.

\section{Approximation of monotone SVFs}\label{sectionMonotone}
Next we obtain several results specific to the approximation of
monotone SVFs by positive sample-based operators. We begin with a
simple condition for the monotonicity preservation by positive
sample-based operators for real-valued functions.
\begin{lemma}\label{MonotoneAverageNumbers}
 Let $\alpha_0,...,\alpha_n \in
[0,1]$, $\beta_0,...,\beta_n \in [0,1]$, $\sum\limits_{i = 0}^n
{\alpha _i} = \sum\limits_{i = 0}^n {\beta_i} = 1 $, such that for
any integer $k$, $0 \leq k \leq n$,
\begin{equation}\label{SequenceDominance}
\sum\limits_{i = k}^n {\alpha _i  \leq } \sum\limits_{i = k}^n
{\beta _i } \ .
\end{equation}
The condition (\ref{SequenceDominance}) is necessary and sufficient
for the inequality,
\begin{equation}\label{SequenceSumNumbers}
\mathop  \otimes \limits_{i = 0}^n \alpha _i r_i  \leq \mathop
\otimes \limits_{i = 0}^n \beta _i r_i \ ,
\end{equation}
to hold for any monotone non-decreasing sequence of numbers
$r_0,...,r_n \in \mathbb{R}$.
\end{lemma}
\proof The necessity follows by considering sequences of the form
$r_0  = r_1  = ... = r_k  < r_{k + 1}  = r_{k + 2} ... = r_n $. The
sufficiency can be obtained by setting,
\begin{equation*}
r_j  = \Delta _0  + \sum\limits_{i = 1}^j {\Delta _i },\ j = 0,...,n
\ ,
\end{equation*}
where $\Delta_0 = r_0$, $\Delta _i = r_i - r_{i-1}>0$, $i =
1,...,n$, and considering the contribution of each $\Delta_i$ in
(\ref{SequenceSumNumbers}). \qed

In case $r_0,...,r_n$ is monotone non-increasing, condition
(\ref{SequenceDominance}) implies by symmetry, that
\begin{equation*}
\mathop  \otimes \limits_{i = 0}^n \alpha _i r_i  \geq \mathop
\otimes \limits_{i = 0}^n \beta _i r_i \ .
\end{equation*}
\begin{corollary}\label{RealMonotonicityPresrving}
Let  $\tilde O_n\left( {F,x} \right)$ be a positive sample-based
operator defined as in (\ref{PositiveOperatorEq}), such that  for
any $x,y \in [0,1]$,$x \leq y$, the weights $\alpha_i = c_{n,i}(x)$,
$\beta_i = c_{n,i}(y)$ satisfy (\ref{SequenceDominance}).  Then
$\tilde O_n$ is monotonicity preserving.
\end{corollary}
Next we show that similar conditions are necessary and sufficient
for the monotonicity preservation by positive sample-based operators
for SVFs. A sequence of sets ${\left\{ {F_i } \right\}_{i \in
\mathbb{Z}} }$ is termed monotone non-decreasing (non-increasing),
if for all $i$, $F_i  \subseteq F_{i + 1} \left( {F_i  \supseteq
F_{i + 1} } \right)$. Monotone non-decreasing (non-increasing) SVFs
are defined in a similar way.

\begin{lemma}\label{MonotoneAverage}
Let $\alpha_0,...,\alpha_n$ , $\beta_0,...,\beta_n $ be as in Lemma
\ref{MonotoneAverageNumbers}, then condition
(\ref{SequenceDominance}) is necessary and sufficient for the
relation,
\begin{equation*}
\mathop  \otimes \limits_{i = 0}^n \alpha _i A_i  \subseteq \mathop
\otimes \limits_{i = 0}^n \beta _i A_i \ ,
\end{equation*}
to hold for any monotone non-decreasing sequence of sets
$A_0,...,A_n \in \mathfrak{J}$.
\end{lemma}
\proof In view of Property \ref{partAverageAllSameProp} of Theorem
\ref{partitionAveargeProps}, the necessity follows by considering
sequences of the form $A_0 = A_1 = ... = A_k  \subset A_{k + 1}  =
A_{k + 2} ... = A_n $. To obtain the sufficiency, assume that
$A_0,...,A_n$ is monotone non-decreasing and let $\Omega _\chi $ be
as in Definition \ref{PartionOfUnion}. We observe that due to the
monotonicity of the sequence $ A_0,...,A_n$, if $\Omega _\chi \ne
\phi$ then necessarily $\chi = \left\{ {k,k + 1,...,n} \right\}$ for
some integer $k$, $0 \leq k \leq n$. Using (\ref{SequenceDominance})
and Definition \ref{SubSetGenFunction} we obtain that for $\Omega
_\chi $ as above,
\begin{equation*}
\left[ {\Omega _\chi  } \right]_{\sum\limits_{i \in \chi } {\alpha
_i } }  \subseteq \left[ {\Omega _\chi  } \right]_{\sum\limits_{i
\in \chi } {\beta _i } } ,
\end{equation*}
which in view of Definition \ref{partitionAvearge} of the partition
average completes the proof of the lemma. \qed

In case $A_0,...,A_n$ is monotone non-increasing, condition
(\ref{SequenceDominance}) implies by symmetry, that
\begin{equation*}
\mathop \otimes \limits_{i = 0}^n \alpha _i A_i  \supseteq \mathop
\otimes \limits_{i = 0}^n \beta _i A_i \ .
\end{equation*}

\begin{corollary}\label{SetMonotonicityPreserving}
Let  $O_n\left( {F,x} \right)$ be a positive sample-based operator
defined as in (\ref{PositiveOperatorEq}), such that  for any $x,y
\in [0,1]$,$x \leq y$, the weights $\alpha_i = c_{n,i}(x)$, $\beta_i
= c_{n,i}(y)$ satisfy (\ref{SequenceDominance}).  Then $O_n$ is
monotonicity preserving.
\end{corollary}
From Corollaries \ref{RealMonotonicityPresrving} and
\ref{SetMonotonicityPreserving} we conclude
\begin{corollary} Let $\widetilde O_n$ and $O_n$ be
defined as in (\ref{RealPositiveOperatorEq}) and
(\ref{PositiveOperatorEq}) respectively . Then $O_n$ is monotonicity
preserving if and only if $\widetilde O_n$ is monotonicity
preserving.
\end{corollary}

To obtain from (\ref{SequenceDominance}) that the Bernstein
approximation operators are monotonicity  preserving, we need to
show that for $0 \leq x_1 \leq x_2  \leq 1$, $0 \leq k < n$,
\begin{equation}\label{CoefEquation}
\sum\limits_{i = k}^n {b\left( {n,x_1 ;i} \right) \leq }
\sum\limits_{i = k}^n {b\left( {n,x_2 ;i} \right)}  \  .
\end{equation}
This can be observed from the properties of the cumulative binomial
distribution \cite{zelen1964probabilitybook},
\begin{equation}\label{BinomialSum}
\sum\limits_{i = 0}^k {b\left( {n,x;i} \right) = \left( {n - k}
\right)\left( {\begin{array}{*{20}c}
   n  \\
   k  \\

 \end{array} } \right)\int\limits_0^{1 - x} {t^{n - k - 1} \left( {1 - t} \right)^k dt} } \ ,
\end{equation}
which is clearly monotone non-increasing in $x$ and thus leads to
(\ref{CoefEquation}).

To continue the discussion, we recall the notion of the speed of a
curve in a metric space (\emph{see e.g.} \cite{burago2001course},
Chapter 2), which indicates the "smoothness" of a set-valued
function. For a real-valued $f$ the speed at a point $x$ is
\begin{equation*}
v_f \left( x \right) = \mathop {\lim }\limits_{\varepsilon \to 0 }
\frac{{\left| {f\left( x \right) - f\left( {x + \varepsilon }
\right)} \right|}} {{\left| \varepsilon  \right|}} \ ,
\end{equation*}
whenever the limit exists. For differentiable $f$, $v_f$ is the
absolute value of the derivative of $f$. The speed of a SVF $F$ is
defined as
\begin{equation*}
v_F \left( x \right) = \mathop {\lim }\limits_{\varepsilon \to 0 }
\frac{{d_\mu  \left( {F\left( x \right),F\left( {x + \varepsilon }
\right)} \right)}} {{\left| \varepsilon  \right|}} \ .
\end{equation*}
Using relation (\ref{IncludedDistance}) and Corollary
\ref{MeasureTransition}, we obtain
\begin{corollary} Let $\widetilde O_n$ and $O_n$ be monotonicity
preserving operators defined as above. Then for a monotone SVF F,
the speed of  $O_n\left( {F, \cdot } \right)$ equals that of
$\widetilde O_n\left( {\mu \left( F \right), \cdot } \right)$.
\end{corollary}
In particular, for the Bernstein set-valued operators applied to a
monotone SVF F, the speed of $B_n\left( {F, \cdot } \right)$ is a
polynomial, since $B_n \left( {\mu \left( F \right),\cdot} \right)$
is a monotone polynomial and therefore its speed is a polynomial
too.

\section{Approximation of multivariate SVFs}\label{sectionMulti}
The results of Section \ref{sectionOperators} can be generalized to
SVFs defined on a compact subset $K$ of $\mathbb{R}^d$. In this case
we adapt to SVFs families of positive sample-based operators of the
form,
\begin{equation}\label{ExRealPositiveOperatorEq}
\widetilde O_n \left( {f,p} \right) = \sum\limits_{i = 0}^{l_n }
{c_{n,i} \left( p \right)f\left( {p_{n,i} } \right)}  \; , \; p \in
K, \; n \in \mathbb{N} \;,
\end{equation}
where $f : K \to \mathbb{R}$, $p_{n,i} \in K, i=0,...,l_n$,
$c_{n,i}(p) \geq 0$ and $\sum\limits_{i = 0}^{l_n } {c_{n,i} \left(
p \right) = 1} $. In analogy with the univariate case, we define $
\delta_n \left(p\right) = \mathop {\min }\limits_{i \in \{ 0,...,l_n
\} } \left|| {p - p_{n,i} } \right|| $, and assume that for any $p
\in K$,
\begin{equation}\label{deltaNCondition}
\mathop {\lim }\limits_{n \to \infty } \delta _n \left( p \right) =
0 \ .
\end{equation}
Notice that (\ref{ExRealPositiveOperatorEq}) includes many well
known families of approximation operators. Some examples are the
approximation by tensor product Bernstein polynomials, tensor
product Schoenberg splines operators (\cite{de2001practical},
Chapter XVII) and multivariate Bernstein polynomials on simplices
(see e.g. \cite{mathe2003asymptotic}). Similarly to the univariate
case, we define for $F : K \to \mathfrak{J}$,
\begin{equation}\label{ExPositiveOperatorEq}
O_n \left( {F,p} \right) = \mathop  \otimes \limits_{i = 0}^{l_n }
c_{n,i} \left( p \right)F\left( {p_{n,i} } \right)  \ .
\end{equation}

With definitions
(\ref{ExRealPositiveOperatorEq})-(\ref{ExPositiveOperatorEq}), the
analogs of Corollaries \ref{MeasureTransition},
\ref{ContinuousApproxExtension} and \ref{HolderApproximationEx} for
multivariate SVFs are easily derived.

As an example of the application of (\ref{ExPositiveOperatorEq}), we
consider the adaptation to SVFs of the piecewise linear
interpolation over triangulations (see, e.g. \cite{farin2002curves},
Chapter 3). We briefly recall that for a collection of points $P =
\left\{ {p_1 ,...,p_l } \right\} \subset \mathbb{R}^2$, the
\emph{triangulation} $\Gamma $ of $P$ is a collection of triangles
such that
\begin{itemize}
\item The vertices of the triangles consist of points in $P$.
\item The interiors of any two triangles do not intersect.
\item If two triangles are not disjoint, then they share either a vertex or an edge.
\item No edge can be added between
points in $P$ without intersecting an edge of one of the triangles
in $\Gamma$.
\end{itemize}

In the notation of (\ref{ExRealPositiveOperatorEq}), let $P_0  =
\left\{ {p_{0,1} ....p_{0.l_0 } } \right\} \subset \mathbb{R}^2$,
and let $K$ be the \emph{convex hull} of the points in $P_0$. Assume
that the sequence $\left\{ {P_n } \right\}_{n \in \mathbb{N}} $
  is nested, $P_0 \subset
P_1 \subset P_2 ...$. Let $\Gamma_n$ be a triangulation of $P_n$,
such that
\begin{equation*}
\mathop {\lim }\limits_{n \to \infty } \Delta _n  = 0 \ ,
\end{equation*}
where $\Delta _n  = \mathop {\max }\limits_{T \in \Gamma _n }
{\text{diam}}\left( T \right)$ and ${\text{diam}}\left( T \right)$
is the diameter of the circumscribed circle of $T$. Note that for
such a sequence of triangulations condition (\ref{deltaNCondition})
is satisfied.

Let the triangle $\Delta = p_{n,{k_1}}, p_{n,{k_2}}, p_{n,{k_3}}$ be
in $\Gamma_n$ and let $p \in \Delta$. The piecewise linear
interpolant $\tilde L_n \left( {f,p} \right)$ is defined as in
(\ref{ExRealPositiveOperatorEq}) with the weights
$c_{n,i}\left(p\right)$ given by
\begin{equation}\label{triangularWeights1}
\begin{array}{*{20}c}
   {c_{n,k_1 } \left( p \right) = \frac{{area(p,p_{n,k_2 } ,p_{n,k_3 } )}}
{{area(p_{n,k_1 } ,p_{n,k_2 } ,p_{n,k_3 } )}},} & {c_{n,k_2 } \left(
p \right) = \frac{{area(p_{n,k_1 } ,p,p_{n,k_3 } )}}
{{area(p_{n,k_1 } ,p_{n,k_2 } ,p_{n,k_3 } )}},}  \\\\
   {c_{n,k_3 } \left( p \right) = \frac{{area(p_{n,k_1 } ,p_{n,k_2 } ,p)}}
{{area(p_{n,k_1 } ,p_{n,k_2 } ,p_{n,k_3 } )}},\;} & {}  \\
\end{array}
\end{equation}
and
\begin{equation}\label{triangularWeights2}
c_{n,i} \left( p \right) = 0, \ i \notin \left\{ {k_1 ,k_2 ,k_3 }
\right\} \ .
\end{equation}
It is easy to verify that for a continuous $f$,
\begin{equation}
\mathop {\lim }\limits_{n \to \infty } \tilde L_n \left( {f,p}
\right) = f\left( p \right) \ ,
\end{equation}
and for $f \in {\text{Lip}}\left( {L,\nu } \right)$,
\begin{equation}
\left| {f\left( p \right) - \tilde L_n \left( {f,p} \right)} \right|
\leq L\Delta _n^\nu  \ .
\end{equation}
Similarly, for a SVF $F$, the piecewise interpolant $L_n \left(
{F,p} \right)$ is defined as in (\ref{ExPositiveOperatorEq}), with
the weights $c_{n,i}(p)$ given by
(\ref{triangularWeights1})-(\ref{triangularWeights2}). Using
Corollaries \ref{ContinuousApproxExtension} and
\ref{HolderApproximationEx} extended to multivariate SVFs, we obtain
that for a continuous $F$
\begin{equation}
\mathop {\lim }\limits_{n \to \infty } L_n \left( {F,p} \right) =
F\left( p \right) \ ,
\end{equation}
and for $F \in {\text{Lip}}\left( {L,\nu } \right)$,
\begin{equation}
\left| {F\left( p \right) -  L_n \left( {F,p} \right)} \right| \leq
2L\Delta _n^\nu  \ .
\end{equation}
Note that in the real-valued case the zero-weighted summands in
(\ref{ExRealPositiveOperatorEq}) do not affect the result, but this
is not so in the set-valued case. More precisely, let $\hat
L_n\left( {F,p} \right) = \mathop  \otimes \limits_{i = 0,c_{n,i}
\ne 0}^{l_n } c_{n,i} \left( p \right)F\left( {p_{n,i} } \right)$.
Then in view of Remark \ref{remarkShort},
\begin{equation*}
L_n \left( {F,p} \right) \ne \hat L_n \left( {F,p} \right) \ .
\end{equation*}
Moreover, while $L_n \left( {F,p} \right)$ is continuous by its
definition, $\hat L_n \left( {F,p} \right)$ is discontinuous in view
of Remark \ref{remarkShort}.

\section{Approximation of functions with values in general metric spaces}\label{sectionGeneral}
Finally, we extend the approximation results for functions with
values in the metric space $\left\{\mathfrak{J}, d_\mu\right\}$ and
the partition average to functions with values in general metric
spaces endowed with an average satisfying certain properties.

Let $\left\{ {X,d_X } \right\}$ be a metric space, and let $
\boxplus$ be an average on elements of $X$ defined for non-negative
weights. Assume that the average $\boxplus$ satisfies the
conditions, that for any $ \Lambda _0 ,...,\Lambda _n \in X$ and
$\alpha _0 ,...,\alpha _n \in [0,1]$, $ \sum\limits_{i = 0}^n
{\alpha _i = 1} $,
\begin{equation}\label{generalMetricCondition}
\mathop  \boxplus \limits_{i = 0}^n \alpha _i \Lambda _i  \in X \ \
\text{and} \ \ d_X\left( {\Lambda _j ,\mathop  \boxplus \limits_{i =
0}^n \alpha_i\Lambda _i } \right) \leq \sum\limits_{i = 0}^n
{\alpha_i d_X \left( {\Lambda _j ,\Lambda _i } \right)} \;, \  j \in
\left\{0,..,n\right\}.
\end{equation}

Let $\widetilde O_n$ be defined by (\ref{ExRealPositiveOperatorEq}),
we define for $G : K \to X$,
\begin{equation*}
O_n \left( {F,x} \right) = \mathop  \boxplus \limits_{i = 0}^n
c_{n,i} \left( x \right)G\left( {x_{n,i} } \right) \; .
\end{equation*}
With these definitions, it is straightforwardly to obtain
approximation results similar to Corollaries
\ref{ContinuousApproxExtension} and \ref{HolderApproximationEx}.

To characterize metrics spaces, in which averages satisfying the
relation (\ref{generalMetricCondition}) can be constructed, we
observe that (\ref{generalMetricCondition}) is equivalent to the
condition that
\begin{equation}\label{ballMetricCondition}
\mathop  \boxplus \limits_{i = 0}^n \alpha _i \Lambda _i  \in
\bigcap\limits_{i = 0}^n {Bl\left( {\Lambda _i ,\sum\limits_{j =
0}^n {\alpha _i d_X \left( {\Lambda _i ,\Lambda _j } \right)} }
\right)} \ ,
\end{equation}
where $Bl\left(\Lambda, r\right)$ is the \emph{metric ball} of
radius $r$ centered at $\Lambda$.

Therefore, we say that a metric space is \emph{strongly convex}, if
for any $ \Lambda _0 ,...,\Lambda _n \in X$ and $\alpha _0
,...,\alpha _n \in [0,1]$, $\sum\limits_{i = 0}^n {\alpha _i }  =
1$,  the set
\begin{equation*}
\Phi \left( {\Lambda _0 ,...,\Lambda _n ;\alpha _0 ,...,\alpha _n }
\right) = \bigcap\limits_{i = 0}^n {Bl\left( {\Lambda _i
,\sum\limits_{j = 0}^n {\alpha _i d_X \left( {\Lambda _i ,\Lambda _j
} \right)} } \right)} \ ,
\end{equation*}
is not empty. Notice that for $n = 1$, the above definition
coincides with the definition of a \emph{convex metric space} in the
sense of Menger (see, e.g., \cite{papadopoulos2005metric}, Chapter
2). In a strongly convex metric space $X$ one can define the average
of any $\Lambda _0 ,...,\Lambda _n \in X$ with the weights $\alpha
_0 ,...,\alpha _n \in [0,1]$, $ \sum\limits_{i = 0}^n {\alpha _i =
1} $ as any element in the set $ \Phi \left( {\Lambda _0
,...,\Lambda _n ;\alpha _0 ,...,\alpha _n } \right)$ .

\bibliographystyle{abbrv}

\end{document}